\documentclass[11pt]{article}
\usepackage{mathtools}

\usepackage[margin=1in]{geometry}
\usepackage{sn-preamble}
\usepackage{latexsym}
\usepackage{amsmath}
\usepackage{amssymb}
\usepackage{amsthm}

\newcommand{\ignore}[1]{}

\newcommand{\Poincare}{Poincar\'e}

\newcommand{\NInf}{\Inf^{-}}

\newcommand{\sens}{\mathrm{sens}}
\newcommand{\nsens}{\mathrm{sens}^{-}}

\title{A Counterexample to a Directed KKL Inequality\vspace{1em}}

\author{
Quentin Dubroff\\[0.15em] \textsl{Rutgers University}
\and Shivam Nadimpalli\\[0.15em] \textsl{Columbia University}
\and Bhargav Narayanan\\[0.15em] \textsl{Rutgers University}
\vspace{1em}
} 

\date{\today}

\begin{document}

\maketitle

\begin{abstract}
We show that the natural directed analogues of the KKL theorem~\cite{KKL:88} and the Eldan--Gross inequality~\cite{Eldan2020} from the analysis of Boolean functions fail to hold. This is in contrast to several other isoperimetric inequalities on the Boolean hypercube (such as the \Poincare~inequality, Margulis's inequality~\cite{Margulis:74} and Talagrand's inequality~\cite{Talagrand:93}) for which directed strengthenings have recently been established.
\end{abstract}

\section{Introduction}
\label{sec:intro}

In this note, we consider isoperimetric inequalities over the Boolean hypercube $\zo^n$. Our notation and terminology follow O'Donnell~\cite{odonnell-book}; in particular, we refer the reader to Chapter~2 of \cite{odonnell-book} for further background.

Recall that given a Boolean function $f\isazofunc$ and an input $x \in \zo^n$, we define the \emph{sensitivity of $f$ at $x$} as 
\[\sens_f(x) := \#\big\{ i : f(x) \neq f(x^{\oplus i})\big\} \qquad \text{where } x^{\oplus i} := (x_1, \ldots, 1-x_i, \ldots, x_n).\]
Two closely related isoperimetric quantities are the \emph{influence of a variable $i\in[n]$ on $f$}, given by 
\[\Inf_i[f] := \Prx_{\bx\sim\zo^n}\sbra{f(\bx)\neq f(\bx^{\oplus i})},\]
and the \emph{total influence of $f$}, given by 
\[\TInf[f] := \sumi \Inf_i[f].\]
It is easy to check that $\TInf[f] = \Ex[\sens_f(\bx)]$, and so the total influence of a function is sometimes also referred to as its \emph{average sensitivity}. 

To set the stage, we recall perhaps the simplest isoperimetric inequality on the Boolean hypercube, the \Poincare~inequality, which says that 
\[\TInf[f] \geq \Var[f].\]
The follow strengthening of the \Poincare~inequality was obtained by Talagrand~\cite{Talagrand:93}, which is known to imply yet another isoperimetric inequality due to Margulis~\cite{Margulis:74}. 

\begin{theorem}[Talagrand's inequality] \label{prop:undirected-talagrand}
	Given a Boolean function $f\isazofunc$, we have 
	\[\Ex_{\bx\sim\zo^n}\sbra{\sqrt{\sens_f(\bx)}}\geq \Omega\pbra{\Var[f]}.\]
\end{theorem}

An alternative (and incomparable) strengthening of the \Poincare~inequality is given by the celebrated Kahn--Kalai--Linial theorem~\cite{KKL:88}. 

\begin{theorem}[KKL inequality] \label{prop:kkl}
	Given a Boolean function $f\isazofunc$, there exists $i\in[n]$ such that 
	\[\Inf_i[f] \geq \Omega\pbra{\Var[f]\cdot\frac{\log n}{n}}.\]
\end{theorem}

Talagrand~\cite{Talagrand1997} conjectured the following common generalization of \Cref{prop:undirected-talagrand,prop:kkl}, which was proved by Eldan~and~Gross~\cite{Eldan2020}.

\begin{theorem}[Eldan--Gross inequality] \label{prop:EG}
	Given a Boolean function $f\isazofunc$, we have 
	\[\Ex_{\bx\sim\zo^n}\sbra{\sqrt{\sens_f(\bx)}}\geq \Omega\pbra{\Var[f]\sqrt{\log\pbra{2 + \frac{e}{\sumi \Inf_i[f]^2}}}}.\]
\end{theorem}

In this note, we will be concerned with directed versions of such results in the Boolean hypercube. Recall that a Boolean function $f\isazofunc$ is said to be \emph{monotone} (resp. \emph{anti-monotone}) if for all $x, y\in\zo^n$, $x\preccurlyeq y$ implies $f(x)\leq f(y)$ (resp. $f(x)\geq f(y)$).\footnote{We write $x\preccurlyeq y$ to mean $x_i\leq y_i$ for all $i\in[n]$.} In connection with the problem of monotonicity testing, Khot, Minzer, and Safra~\cite{kms} obtained a {``directed''} analogue of \Cref{prop:undirected-talagrand}. We write
\[\nsens_f(x) := \#\big\{i : f(x)> f(x^{\oplus i}) \text{ and } x \preccurlyeq x^{\oplus i}\big\}\] 
for the \emph{negative sensitivity of $f$ at $x$}, and write 
\[
\eps(f) := \min_{g \text{ monotone}} \dist(f,g) \qquad\text{where}\qquad \dist(f,g) := \Prx_{\bx\sim\zo^n}\sbra{f(\bx)\neq g(\bx)}
\]
for the \emph{distance to monotonicity of $f$}.

\begin{theorem}[Theorem 1.6~of~\cite{kms}] \label{prop:dirTal}
	Given a Boolean function $f\isazofunc$, we have 
	\[\Ex_{\bx\sim\zo^n}\sbra{\sqrt{\nsens_f(\bx)}}\geq \Omega\pbra{\eps(f)}.
	\footnote{The original result due to Khot, Minzer, and Safra had additional logarithmic factors in $n$ and $1/\epsilon(f)$, but this was improved by Pallavoor, Raskhodnikova, and Waingarten~\cite{Pallavoor2022}.}
	\]
\end{theorem}

Indeed, prior results on monotonicity testing due to Goldreich et al.~\cite{GGLRS} and Chakrabarty and Seshahdri~\cite{Chakrabarty2016} can be viewed as directed analogues of the \Poincare~inequality and Margulis's inequality~\cite{Margulis:74} respectively. Finally, a directed analogue of an inequality due to Pisier~\cite{Pisier:86} was obtained by Canonne~et al.~\cite{Canonne2021}.

Although the directed analogues are known to imply their undirected counterparts (cf. Section 9.4 of~\cite{kms}), their proofs bear little resemblance to the proofs in the undirected setting (with the exception of the directed Pisier inequality) and are usually much more involved. 


These results suggest an informal analogy between the undirected and the directed cube, with isoperimetric quantities being replaced with their directed counterparts and $\Var[f]$ being replaced with $\eps(f)$ in the latter. Writing
\[\NInf_i[f]:= \#\big\{x : f(x) > f(x^{\oplus i}) \text{ and } x \preccurlyeq x^{\oplus i}\big\}\cdot\frac{1}{2^{n-1}}\]
for the \emph{negative influence of $i$ on $f$}, we have the following natural directed analog of \Cref{prop:kkl}.

\begin{conjecture}[Directed KKL inequality] \label{conj:dirKKL}
	Given a Boolean function $f\isazofunc$, there exists $i\in[n]$ such that 
	\[\NInf_i[f] \geq \Omega\pbra{\eps(f)\cdot\frac{\log n}{n}}.\]
\end{conjecture}

\Cref{conj:dirKKL}, as well as a Fourier analytic reformulation thereof, appears to have been raised by Subhash Khot at the April~2016 Simons Meeting on Algorithms and Geometry~\cite{Holden}. Our aim in this short note is to show that \Cref{conj:dirKKL} fails to hold. 

\begin{theorem}\label{th:main}
There is a function $f:\zo^{2n}\to\zo$ with 
\begin{enumerate}
	\item $\NInf_i[f] = 0$ for all $i\in[n]$,
	\item $\NInf_i[f] = O(1/n)$ for all $i\in[2n] \setminus [n]$, and
	\item $\eps(f) = \Omega(1)$,
\end{enumerate}
\end{theorem}
We note that this further rules out a natural directed analog of \Cref{prop:EG} (which would imply \Cref{conj:dirKKL}). The construction establishing \Cref{th:main} follows in the next section.

\begin{remark}
After a draft of this paper was circulated, it was brought to our attention that Minzer and Khot~\cite{DorEmail} have independently discovered a similar construction to the one establishing \Cref{th:main}.
\end{remark}



\section{A Counterexample to Directed KKL}
\label{sec:counterexample}

We view $\zo^{2n}$ as $\zo^{n}\times\zo^n$ and construct a function $f:\zo^n\times\zo^n\to\zo$ with 
\begin{enumerate}
	\item $\NInf_i[f] = 0$ for all $i\in[n]$,
	\item $\NInf_i[f] = O(1/n)$ for all $i\in[2n] \setminus [n]$, and
	\item $\eps(f) = \Omega(1)$,
\end{enumerate}
thereby refuting \Cref{conj:dirKKL}. 

\begin{proof}[Proof of \Cref{th:main}] Let $T_1, \ldots, T_n \in {[n]\choose {\log n}}$ be drawn independently and uniformly at random. Set 
\[f(x, y) := \bigvee_{i=1}^n \pbra{\pbra{\bigwedge_{j\in T_i} x_j} \wedge (1-y_i)}.\]
We note that this function is closely related to the well-known ``Tribes'' function due to Ben-Or and Linial~\cite{BenOrLinial:85short}.

It is clear that $f$ is monotone in the first $n$ coordinates and anti-monotone in the last $n$ coordinates; consequently $\NInf_i[f] = 0$ for all $i\in [n]$. A coordinate $i \in [2n]\setminus[n]$ is relevant only on $x\in\zo^n$ for which $\bigwedge_{j\in T_i} x_j = 1$; as $|T_i| = \log n$, this set has measure at most 
\[\frac{2^{n-\log n}}{2^n} = \frac{1}{n}.\]
It follows that $\NInf_i[f] = O(1/n)$ for all $i\in[2n]$. 

Before turning to the third item above, we recall the following fact from \cite{kms} without proof.

\begin{lemma}[Lemma~3.11 of \cite{kms}]
	For $f:\zo^n\times\zo^n\to\zo$ such that $f$ is monotone in the first $n$ coordinates and anti-monotone in the last $n$ coordinates, we have 
	\[\eps(f) = \Theta\pbra{\Ex_{\bx\sim\zo^n}\sbra{\Varx_{\by\sim\zo^n}\sbra{f(\bx,\by)}}}. \eqno
	\]
\end{lemma}

Suppose, for convenience, that $x\in\zo^n$ is such that $\bigwedge_{j\in T_i} x_j = 1$ for exactly one $i\in[n]$. Then the restricted function $f(x,\cdot):\zo^n\to\zo$ is simply the anti-dictatorship $(1-y_i)$, and has $\Var\sbra{f(x,\cdot)} = \Omega(1)$. We will be done if we can show that this happens for $\Omega(1)$ fraction of $x\in\zo^n$. As before, for fixed $i\in [n]$ we have 
\[\Prx_{\bx\sim\zo^n}\sbra{\bigwedge_{j\in T_i} \bx_j} = \frac{1}{n}
\qquad\text{and so}\qquad
\Ex_{\bx\sim\zo^n}\sbra{\#\bigg\{i : \bigwedge_{j\in T_i} \bx_j = 1\bigg\}} = 1.
\]
By Markov's inequality, we thus have 
\[\Prx_{\bx\sim\zo^n}\sbra{\#\bigg\{i : \bigwedge_{j\in T_i} \bx_j = 1\bigg\} \geq 2} \leq \frac{1}{2}.\]
We also have 
\[\Prx_{\bx\sim\zo^n}\sbra{\#\bigg\{i : \bigwedge_{j\in T_i} \bx_j = 1\bigg\} = 0} \approx \pbra{1 - \frac{1}{n}}^n \approx \frac{1}{e},\]
and so the desired event happens with constant probability, and we are done.
\end{proof}

\bibliography{allrefs.bib}

\newcommand{\etalchar}[1]{$^{#1}$}
\begin{thebibliography}{CCK{\etalchar{+}}21}

\bibitem[BOL85]{BenOrLinial:85short}
M.~Ben-Or and N.~Linial.
\newblock Collective coin flipping.
\newblock In {\em Proc. 26th Annual Symposium on Foundations of Computer
  Science (FOCS)}, pages 408--416, 1985.

\bibitem[CCK{\etalchar{+}}21]{Canonne2021}
Cl{\'{e}}ment~L. Canonne, Xi~Chen, Gautam Kamath, Amit Levi, and Erik
  Waingarten.
\newblock Random restrictions of high dimensional distributions and uniformity
  testing with subcube conditioning.
\newblock In D{\'{a}}niel Marx, editor, {\em Proceedings of the 2021 {ACM-SIAM}
  Symposium on Discrete Algorithms, {SODA} 2021, Virtual Conference, January 10
  - 13, 2021}, pages 321--336. {SIAM}, 2021.

\bibitem[CS16]{Chakrabarty2016}
Deeparnab Chakrabarty and C.~Seshadhri.
\newblock An $o(n)$ monotonicity tester for boolean functions over the
  hypercube.
\newblock {\em {SIAM} J. Comput.}, 45(2):461--472, 2016.

\bibitem[EG20]{Eldan2020}
Ronen Eldan and Renan Gross.
\newblock Concentration on the boolean hypercube via pathwise stochastic
  analysis.
\newblock In Konstantin Makarychev, Yury Makarychev, Madhur Tulsiani, Gautam
  Kamath, and Julia Chuzhoy, editors, {\em Proccedings of the 52nd Annual {ACM}
  {SIGACT} Symposium on Theory of Computing, {STOC} 2020, Chicago, IL, USA,
  June 22-26, 2020}, pages 208--221. {ACM}, 2020.

\bibitem[GGL{\etalchar{+}}00]{GGLRS}
O.~Goldreich, S.~Goldwasser, E.~Lehman, D.~Ron, and A.~Samordinsky.
\newblock Testing monotonicity.
\newblock {\em Combinatorica}, 20(3):301--337, 2000.

\bibitem[KKL88]{KKL:88}
J.~Kahn, G.~Kalai, and N.~Linial.
\newblock The influence of variables on boolean functions.
\newblock In {\em Proc. 29th Annual Symposium on Foundations of Computer
  Science (FOCS)}, pages 68--80, 1988.

\bibitem[KMS15]{kms}
Subhash Khot, Dor Minzer, and Muli Safra.
\newblock On monotonicity testing and boolean isoperimetric type theorems.
\newblock In {\em 2015 IEEE 56th Annual Symposium on Foundations of Computer
  Science}, pages 52--58, 2015.

\bibitem[Lee22]{Holden}
Holden Lee.
\newblock Notes on {S}imons {A}lgorithms~and~{G}eometry {M}eetings.
\newblock 2022.
\newblock
  \href{https://www.dropbox.com/s/sc8komdxy4ofs8j/simons.pdf?dl=0}{Link}.

\bibitem[Mar74]{Margulis:74}
G.~Margulis.
\newblock Probabilistic characteristics of graphs with large connectivity.
\newblock {\em Prob.\ Peredachi Inform.}, 10:101--108, 1974.

\bibitem[Min22]{DorEmail}
Dor Minzer.
\newblock Personal communication.
\newblock 2022.

\bibitem[O'D14]{odonnell-book}
R.~O'Donnell.
\newblock {\em Analysis of Boolean Functions}.
\newblock Cambridge University Press, 2014.

\bibitem[Pis86]{Pisier:86}
G.~Pisier.
\newblock Probabilistic methods in the geometry of {B}anach spaces.
\newblock In {\em Lecture notes in Math.}, pages 167--241. Springer, 1986.

\bibitem[PRW22]{Pallavoor2022}
Ramesh Krishnan~S. Pallavoor, Sofya Raskhodnikova, and Erik Waingarten.
\newblock Approximating the distance to monotonicity of boolean functions.
\newblock {\em Random Struct. Algorithms}, 60(2):233--260, 2022.

\bibitem[Tal93]{Talagrand:93}
M.~Talagrand.
\newblock Isoperimetry, logarithmic {S}obolev inequalities on the discrete cube
  and {M}argulis' graph connectivity theorem.
\newblock {\em GAFA}, 3(3):298--314, 1993.

\bibitem[Tal97]{Talagrand1997}
M.~Talagrand.
\newblock On boundaries and influences.
\newblock {\em Combinatorica}, 17(2):275--285, 1997.

\end{thebibliography}
\bibliographystyle{alpha}

\appendix

\end{document}